\documentclass[12pt,a4paper,reqno]{amsart}
\usepackage{amsfonts, amsthm, amssymb, amsmath}

\theoremstyle{plain}
\newtheorem*{theorem}{Theorem}
\newtheorem*{corollary}{Corollary}
\newtheorem{lemma}{Lemma}
\newtheorem{extralemma}{Lemma~A.\ignorespaces}

\theoremstyle{definition}
\newtheorem{definition}{Definition}
\newtheorem{remark}{Remark}
\newtheorem{extradefinition}{Definition~A.\ignorespaces}
\newtheorem*{scc}{The condition $C'(\lambda)$}

\begin{document}
\title[Consistent Solution of Markov's Problem]
{Consistent Solution of Markov's Problem about Algebraic Sets}

\author{Ol'ga~V.~Sipacheva}
\thanks{This work was financially supported by the Russian Foundation
for Basic Research,
project no.~06-01-00764}

\address
{Department of Mechanics and Mathematics\\
Moscow State University\\
Leninskie Gory\\
Moscow, 119992 Russia}


\email
{o-sipa@yandex.ru}

\begin{abstract}
It is proved that the continuum hypothesis implies the existence of
a group $M$ containing a nonalgebraic unconditionally closed set, i.e.,
a set which is closed in any Hausdorff group topology on $M$ but is not
an intersection of finite unions of solution sets of equations in $M$.
\end{abstract}

\maketitle

\begin{definition}[Markov~\cite{Markov1945}]
A subset $A$ of a group $G$ is said to be \emph{unconditionally
closed} in $G$ if it is closed in any Hausdorff group topology on $G$.
\end{definition}

Clearly, all solution sets of equations in $G$, as well as
their finite unions and arbitrary intersections, are
unconditionally closed. Such sets are called algebraic.
The precise definition is  as follows.

\begin{definition}[Markov~\cite{Markov1945}]
A subset $A$ of a group $G$ with identity element $1$ is said to be
\emph{elementary algebraic} in $G$ if there exists a word $w= w(x)$ in the
alphabet $G\cup \{x^{\pm1}\}$ ($x$ is a variable) such that
$$
A =\{x\in G: w(x) = 1\}.
$$
Finite unions of elementary algebraic sets are called \emph{additively
algebraic} sets.  An arbitrary
intersection of additively algebraic sets is said to be \emph{algebraic}.
Thus, the algebraic sets in $G$ are the solution sets of arbitrary
conjunctions of finite disjunctions of equations.
\end{definition}

In his 1945 paper~\cite{Markov1945}, A.~A.~Markov showed that any
algebraic set is unconditionally closed and
posed the problem of whether the converse is true.
In~\cite{Markov1946} (see also~\cite{Markov1944}), he solved this
problem for countable groups by proving that any unconditionally
closed set in a countable group is algebraic.
The answer is also positive for subgroups of direct products of
countable groups~\cite{new}.

Markov's problem is closely related to the topologizability of
groups. Recall that a group is said to be
\emph{topologizable} if it admits a nondiscrete Hausdorff group topology.
Groups that are not topologizable are called \emph{nontopologizable}.
The problem of the existence of a nontopologizable group was posed by
Markov in the same 1945 paper~\cite{Markov1945}; it was solved under CH
by Shelah in~1976 (published in~1980~\cite{Shelah}). The first
ZFC example was given by Hesse in 1979~\cite{Hesse}; a year later,
Ol'shanskii constructed a countable nontopologizable group in
ZFC~\cite{Ol'shanskii}. More recent results can be found in~\cite{nontop}.

In this paper, we prove the following theorem.

\begin{theorem}
Under CH, there exists a group containing a
nonalgebraic unconditionally closed set.
\end{theorem}

\begin{proof}
Such a group is the nontopologizable group $M$ constructed by
Shelah~\cite{Shelah}. It has many remarkable properties. What we need
is
$$
M= \bigcup_{\alpha\in\omega_1}M_\alpha,
$$
where each $M_\alpha$ is a countable subgroup of
$M_{\beta}$ for any $\beta\ni\alpha$ and all of the
$M_\alpha$ (except possibly $M_1$)
are increasing unions of topologizable
subgroups. The following general observation shows that this is
sufficient for $M$ to have a nonalgebraic unconditionally closed
subset.

\begin{lemma}\label{lemma1}
If $G$ is a nontopologizable group
and any finite subset of $G$ is contained
in a topologizable subgroup of $G$,
then $G\setminus \{1\}$ is a nonalgebraic unconditionally
closed subset of $G$.
\end{lemma}

\begin{proof}
Since $G$ admits no
nondiscrete Hausdorff group topology, the set $A=G\setminus\{1\}$ is
unconditionally closed in $G$. Suppose that it is algebraic. Then
$A=\bigcap_{\gamma \in \Gamma}A_\gamma$, where $\Gamma$ is an
arbitrary index set and each $A_\gamma$ is an additively algebraic
set in $G$. All of the sets $A_\gamma$ must contain $G\setminus \{1\}$;
hence each of them must coincide with $G$ or $G\setminus\{1\}=A$. Clearly,
some of these sets does not contain $1$; thus, $A=A_\gamma$ for
some $\gamma$. This means that $A=\bigcup_{i\le
k}A_i$, where $k\in \omega$ and each $A_i$ is an elementary
algebraic set. This means that there exist words $w_1(x), \dots,
w_k(x)$ in the alphabet $G\cup \{x^{\pm1}\}$ such that
$$%
A_i=\{x\in G:  w_i(x) = 1\}
$$
for $i\le k$. Since the number of letters in each
word is finite, we can
find a topologizable subgroup $H\subset G$ such that all of
the $w_i(x)$ are words in the alphabet $H\cup \{x^{\pm1}\}$.
Thus, the $A_i\cap H$ are
elementary algebraic sets in $H$, and $A\cap
H=H\setminus\{1\}$ is an algebraic (and,
therefore, unconditionally closed) set in $H$, which
contradicts the topologizability of $H$.
\end{proof}

\begin{remark}
Combining Lemma~\ref{lemma1} with the theorem of Markov about
unconditionally closed sets in countable groups, we see that any countable
group which is an increasing union of topologizable subgroups
is topologizable. In particular, all of the groups $M_\alpha$, except
possibly $M_1$, are topologizable, and the group $M$ is uncountable.
\end{remark}

This essentially completes the proof of the theorem. It only remains
to verify that $M$ has sufficiently many topologizable
subgroups.\footnote{It is mentioned in~\cite{Shelah} without proof that
all countable subgroups of $M$ are topologizable. This is not so unless
special care is taken; at least, the group $M_0$, which is
the basis of the inductive construction of $M$, must be
topologizable.}
This requires knowledge
of the structure of the groups $M_\alpha$. Below, we reproduce
(or, to be more precise, reconstruct) the part of Shelah's proof
containing the construction of these groups, which is far from being
overloaded with details, in contrast to misprints and lacunae.
The description of Shelah's group suggested below
slightly differs from that given
in~\cite{Shelah}, but the essence is the same.
The proof uses the notions of a malnormal subgroup and
good fellows over a subgroup. Recall that a subgroup $H$ in a group $G$
is said to be \emph{malnormal} if $g^{-1}Hg\cap
H=\{1\}$ for any $g\in G\setminus H$. Shelah calls
two elements $x$ and $y$ of a group $G$ \emph{good fellows}\footnote{In the
definition of good fellows given by Shelah in~\cite[p.~377]{Shelah},
``$G-H$'' should read ``$H-G$''.} over a
subgroup $H\subset G$ if $x, y\in G\setminus H$ and the double cosets
$Hx^\varepsilon H$ and $Hy^\delta H$ are disjoint for $\varepsilon,
\delta= \pm1$, i.e., $x\notin Hy^{\pm 1}H$. Other
algebraic notions, constructions, and facts used in the proof
are collected in the appendix; the very basic definitions
can be found in~\cite{Schupp}.

The groups $M_\alpha$ are constructed by induction as follows. Let
$$
\{S_\gamma: \gamma\in\omega_1\}
$$
be the family of all infinite countable
subsets of $\omega\times \omega_1$
enumerated in such a way that $S_0=\omega\times\{0\}$
(recall that we have assumed $\mathfrak c= \omega_1$).
Let $M_0$ be the trivial group.
For $M_1$ we
take an arbitrary non--finitely generated countable group
and identify it (as a set)
with $\omega\times \{0\}$.  Suppose that $\alpha\in\omega_1$
and $M_\alpha$ is already constructed.  We identify it with
$\omega\times\alpha$ (each ordinal is considered as the set of all
smaller ordinals).  Let us construct $M_{\alpha+1}$.

Consider the set
$$
T_\alpha=\omega^2\times \alpha\times M_\alpha
$$
of all triples $((i, j), \gamma, h)$, where
$i, j\in \omega$, $\gamma \in \alpha$, and $h\in M_\alpha$. This set is
countable. Let us enumerate it:
$$
T_\alpha=\{((i_n, j_n), \gamma_n, h_n):n\in \omega\};
$$
we require that $i_n\le n$ for any $n$. (Certainly, each of $i$,
$j$, $\gamma$, and $h$ occurs in $T_\alpha$ infinitely many times.)
First, we construct increasing sequences of
countable groups $H^\alpha_n$ and $L^\alpha_n$
such that
\begin{enumerate}
\item[(i)]
each $H^\alpha_n$ is a finitely generated subgroup of $M_\alpha$, and
$\bigcup H^\alpha_n = M_\alpha$;
\item[(ii)]
each $H^\alpha_n$ is a subgroup of $L^\alpha_n$,
each $L^\alpha_n$ is a subgroup of $L^\alpha_{n+1}$,
and
$L^\alpha_n \cap M_\alpha=H^\alpha_n$;
\item[(iii)]
the set $L^\alpha_n\setminus M_\alpha$ is infinite, and its elements are
indexed by pairs of integers:
$$
L^\alpha_n\setminus
M_\alpha=L^\alpha_n\setminus
H_n^\alpha=\{a_{(i,j)}: i\le n,\ j\in \omega\};
$$
\item[(iv)]
if $S_{\gamma_n}\subset M_\alpha$ and $S_{\gamma_n}$ is contained in no
finitely generated subgroup of $M_\alpha$, then
$$
h_n\in
\bigl((S_{\gamma_n}\cap H_{n+1}^\alpha)\cup \{a_{(i_n,j_n)}\}\bigr)^{10000}
\subset L^\alpha_{n+1}
$$
(this is the usual power of a set in the group $L^\alpha_{n+1}$);
\item[(v)]
$H^\alpha_n$ is a malnormal subgroup of $L^\alpha_n$, i.e.,
$a^{-1}H^\alpha_na\cap
H^\alpha_n=\{1\}$ for any $a\in L^\alpha_n\setminus H^\alpha_n$.
\end{enumerate}

The groups $L^\alpha_n$ and $H^\alpha_n$ are defined by induction. We set
$H^\alpha_0=\{1\}$ and let $L^\alpha_0$ be an infinite cyclic group
having trivial intersection with $M_\alpha$. We somehow enumerate the
elements of $L^\alpha_0\setminus \{1\}$ by pairs from $\{0\}\times \omega$:
$$
L^\alpha_0\setminus \{1\}=\{a_{(0, j)}:j\in \omega\}.
$$
Suppose that $H^\alpha_n$ and $L^\alpha_n$ are constructed and
$$
L^\alpha_n\setminus M_\alpha=\{a_{(i, j)}: i\le n,\ j\in \omega\}.
$$
Let us construct $H^\alpha_{n+1}$ and $L^\alpha_{n+1}$. Recall that we have
enumerated all infinite countable subsets of
$\omega\times \omega_1$ at the very
beginning of the construction and that $M_\alpha$ is identified with
$\omega\times\alpha$. If the set $S_{\gamma_n}$ (the
$\gamma_n$ is from the enumeration of the set $T_\alpha$ of triples)
is not contained  in $M_\alpha$ or is contained in a finitely generated
subgroup of $M_\alpha$, then we set
$H^\alpha_{n+1} = \langle H^\alpha_n, h_n\rangle$
(this is the subgroup generated by $H^\alpha_n$ and $h_n$ in
$M_\alpha$; it is finitely generated by the induction hypothesis)
and $L^\alpha_{n+1}=L^\alpha_n \mathbin{*_{H^\alpha_n}} H^\alpha_{n+1}$
(this is the free product of $L^\alpha_n$ and $H^\alpha_{n+1}$
with amalgamation over $H^\alpha_n$; see the appendix).
Otherwise, i.e., if $S_{\gamma_n}$ is contained  in $M_\alpha$
and is not contained in any finitely generated subgroup of
$M_\alpha$, then there exist $x,y\in S_{\gamma_n}\setminus H^\alpha_n$
such that $x\notin H^\alpha_n y^{\pm1}H^\alpha_n\cup h_n H^\alpha_n$ in
$M_\alpha$ (in particular, $x$ and $y$ are good fellows over
$H^\alpha_n$). The proof is similar to that
of Fact~2.2(ii) from~\cite{Shelah}: if
any element of $S_{\gamma_n}\setminus H^\alpha_n$ would belong to
$H^\alpha_n z H^\alpha_n \cup H^\alpha_n z^{-1} H^\alpha_n\cup h_n H^\alpha_n$,
where $z$ is an arbitrary element of $S_{\gamma_n}\setminus H^\alpha_n$,
then $S_{\gamma_n}$ would be contained in the set
$H^\alpha_n z H^\alpha_n \cup H^\alpha_n z^{-1} H^\alpha_n\cup
h_n H^\alpha_n\cup H^\alpha_n$,
which is in turn
contained in a finitely generated subgroup, because $H^\alpha_n$ is finitely
generated (by the induction hypothesis).
In this case, we set
$$%
H^\alpha_{n+1} = \langle H^\alpha_n, x,y, h_n\rangle
$$%
(this subgroup is finitely generated).
Recall that $T_\alpha$ is indexed in such a way that
$i_n\le n$, so the element $a_{(i_n, j_n)}\in L_{i_n}$ is already
defined, and that $H_n^\alpha$ is malnormal in $L_n^\alpha$
by the induction hypothesis. Moreover, by construction,
$h_n^{-1}x\in H_{n+1}^\alpha\setminus H_n^\alpha$. We set $h=h^{-1}_n x$
and consider the word
\begin{multline*}
r_0=h a_{(i_n, j_n)}y a_{(i_n,
j_n)} x a_{(i_n, j_n)} (ya_{(i_n, j_n)})^2 xa_{(i_n, j_n)} (ya_{(i_n,
j_n)})^3\\ \dots xa_{(i_n, j_n)} (ya_{(i_n, j_n)})^{80}\in
L^\alpha_n \mathbin{*_{H^\alpha_n}} H^\alpha_{n+1}.
\end{multline*}
Let $N$ be the normal subgroup generated by this word in
$L^\alpha_n \mathbin{*_{H^\alpha_n}} H^\alpha_{n+1}$. We set
\begin{multline*}
L^\alpha_{n+1}= (L^\alpha_n \mathbin{*_{H^\alpha_n}}
H^\alpha_{n+1})/N
=\langle L^\alpha_n \mathbin{*_{H^\alpha_n}} H^\alpha_{n+1}\mid\\
h_n=x a_{(i_n, j_n)}y
a_{(i_n, j_n)} x a_{(i_n, j_n)} (ya_{(i_n, j_n)})^2  \dots
xa_{(i_n, j_n)} (ya_{(i_n, j_n)})^{80} \rangle
\end{multline*}
(this is the amalgamated  free product of
$L^\alpha_n$ and $H^\alpha_{n+1}$ with one defining relation $r_0=1$).
According to Lemma~A.\ref{scc} and the paragraph after
this lemma in the appendix, the groups $L^\alpha_n$ and
$H^\alpha_{n+1}$ are naturally embedded in
$L^\alpha_{n+1}$ as subgroups, and hence $L^\alpha_n\cap
H^\alpha_{n+1} = H^\alpha_n$; moreover, by Lemma~A.\ref{malnormal}
from the appendix, $H^\alpha_{n+1}$ is
malnormal in $L^\alpha_{n+1}$.  Let us somehow enumerate the elements
of $L^\alpha_{n+1}\setminus (L^\alpha_n\cup M_\alpha)$ by the
elements of $\{n+1\}\times \omega$.

The construction of the groups $H^\alpha_n$ and $L^\alpha_n$ is
completed. The $H^\alpha_n$ satisfy condition (i) because $h_n\in
H^\alpha_n$ for every $n$ and $\{h_n: n\in \omega\}=M_\alpha$ by the
definition of $T_\alpha$. The remaining conditions (ii)--(v) hold by
construction (10000 is taken as an upper bound for the length of the
word $r_0$).

We set $M_{\alpha+1}=\bigcup L^\alpha_n$.

Finally, we define $M_\beta=\bigcup_{\alpha\in\beta} M_\alpha$ for
limit $\beta$ and set $M=\bigcup_{\alpha\in\omega_1} M_\alpha$.

We have constructed the required group $M$. As mentioned, it has many
remarkable properties. In particular, each $M_\alpha$ is a malnormal
subgroup of $M$ (i.e.,
$h^{-1}M_\alpha h \cap M_\alpha = \{1\}$ for any
$h\in M\setminus M_\alpha$) and $S^{10000}=M$ for any uncountable
$S\subset M$ (see Lemma~\ref{add} below).
This immediately implies that $M$ admits no
nondiscrete Hausdorff group topology.  Indeed, suppose that such a
topology exists. Take an arbitrary neighborhood $U$ of the identity
element and consider a neighborhood $V$ for which $V^{10000}\subset
U$. If $V$ is countable, then it is contained in some $M_\alpha$ and,
since $M_\alpha$ is malnormal in $M$, $h^{-1} Vh\cap V=\{1\}$ for any
$h \in M\setminus M_{\alpha}$; thus, $\{1\}$ is an open set,
which contradicts the nondiscreteness of the topology.
Hence $V$ must be uncountable, and $M=V^{10000}\subset U$.

\begin{lemma}\label{add}
Each $M_\alpha$ is a malnormal subgroup of $M$ and
$S^{10000}=M$ for any uncountable
$S\subset M$.
\end{lemma}

The malnormality of $M_\alpha$ in $M$ easily follows from the
construction. Indeed,  it
is sufficient to show that $M_\alpha$ is malnormal in $M_{\alpha+1}$
for each $\alpha$. If $h\in M_{\alpha+1}\setminus M_\alpha$ and
$h^{-1} M_\alpha h \cap M_\alpha\ne \{1\}$,
then there exist $k,l,m\in \omega$ and
$a,b\in M_\alpha$ such that $h\in L^\alpha_k\setminus M_\alpha$, $a\in
H^\alpha_l$, $b\in H^\alpha_m\setminus \{1\}$, and $h^{-1}ah=b$. For
$n=\max\{k,l,m\}$, we have $h\in L^\alpha_n\setminus M_\alpha
= L^\alpha_n\setminus H_n^\alpha$, $a\in
H^\alpha_n$, and $b\in H^\alpha_n\setminus \{1\}$; thus,
$h^{-1}H^\alpha_nh\cap H^\alpha_n\ne \{1\}$, which contradicts~(v).

Let us prove that $S^{10000}=M$ for any uncountable $S$.
First, note that if $S\subset M$ is
uncountable, then there exists a $\beta$ such that
$S\cap M_\beta$ is contained in no finitely
generated subgroup of $M_\beta$. Indeed, take an increasing sequence
of countable ordinals $\beta_k$ such that $S\cap M_{\beta_0}\ne
\varnothing$ and $S\cap M_{\beta_{k+1}}\setminus M_{\beta_{k}}\ne
\varnothing$ for any $k$. Let $\beta=\sup\{\beta_k\}_{k=0}^\infty$.
By definition, $M_\beta=\bigcup_{\lambda\in \beta} M_\alpha$. Any
subgroup of $M_\beta$ generated by finitely many elements $g_1,
\dots, g_n$ is contained in $M_\alpha$ for some $\alpha<\beta$ and,
therefore, in $M_{\beta_k}$ for some $k$. Thus, $S$ is not contained
in any finitely generated subgroup of $M_\beta$.
According to Fact~2.8 in~\cite{Shelah}, $S$
is not contained in any finitely generated subgroup of
$M_\alpha$ for any $\alpha\ni\beta$.
We have $S\cap M_\beta=S_\gamma$ for some
$\gamma$.  Take any $h\in M$ (then $h\in M_\delta$ for some
$\delta$). Since $S$ is uncountable, there exists an
$\alpha\ni\max\{\beta, \gamma, \delta\}$ such that $S\cap
(M_{\alpha+1}\setminus M_\alpha)\ne \varnothing$. Let $a\in
S\cap (M_{\alpha+1}\setminus M_\alpha)$. Then
$a\in L_k^\alpha\setminus M_\alpha$ for some $k$ and, by~(iii),
$a=a_{(i,j)}$ for some
$(i,j) \in \omega^2$ ($i\le k$).  We have
$( (i,j), \gamma, h) \in T_\alpha$,
i.e., $( (i,j), \gamma, h)=
( (i_n,j_n), \gamma_n, h_n)$ for some $n$; in particular, $a_{(i, j)}
=a_{(i_n, j_n)}$, $S_{\gamma}=S_{\gamma_n}$, and $h= h_n$.
The set $S_{\gamma_n}=S_\gamma=S\cap M_\beta$ is contained in
$M_\alpha\supset M_\beta$ but not in a finitely generated subgroup of
$M_\alpha$; hence, by the construction of $L_{n+1}^\alpha$, there
exist $x, y\in S_{\gamma_n}\subset S$ such that
$h_n=
x a_{(i_n, j_n)}y
a_{(i_n, j_n)} x a_{(i_n, j_n)} (ya_{(i_n, j_n)})^2  \dots
xa_{(i_n, j_n)} (ya_{(i_n, j_n)})^{80}$ in
$L^\alpha_{n+1}$ (and in $M$). Thus, $h=h_n$ is a product of
length less than 10000 of elements of $S$.
\end{proof}

It remains to prove that $M$ has sufficiently many
topologizable subgroups. It suffices to show
that, for any $\alpha\in \omega_1\setminus\{0\}$ and $k\in
\omega$, there exists an $n\ge k$ such that the group $L^\alpha_n$ is
topologizable. This is implied by Lemma~A.\ref{lemma2}
from the appendix. Indeed,
note that, for any $\alpha \in \omega_1\setminus\{0\}$
and $k\in \omega$, there
exists an $n\ge k$ such that the group
$H^\alpha_{n+1}$ contains a pair of goods fellows over $H^\alpha_n$,
because, according to Fact~2.8 in~\cite{Shelah}, any set $S$
not contained in a finitely generated subgroup of some $M_\alpha$
is not contained in any finitely generated subgroup of
$M_\beta$ for $\beta>\alpha$.  The group $M_1$ is not finitely
generated; therefore, it is not contained in a finitely generated
subgroup of any of the groups $M_\alpha$. On the other hand,
$M_1=\omega\times \{0\} =S_0$. Each ordinal $\gamma\in\alpha$ occurs
in infinitely many triples from $T_\alpha$; take a triple containing
$\gamma=0$ and having number $n(k)\ge k$ in the enumeration of
$T_\alpha$.  By construction, the group $H^\alpha_{n(k)+1}$ is
generated by $H^\alpha_{n(k)}$, some element $t$ of $M_\alpha$, and
a pair of goods fellows $x,y$ over $H^\alpha_{n(k)}$, for which
$t^{-1}x=h\in H^\alpha_{n(k)+1}\setminus H^\alpha_{n(k)}$; moreover,
there exists an $a\in L^\alpha_{n(k)}$ such that
$L^\alpha_{n(k)+1} =\langle L^\alpha_{n(k)}
\mathbin{*_{H^\alpha_{n(k)}}} H^\alpha_{n(k)+1}\mid
r_0=1\rangle$, where $r_0$ is the same word as in Lemma~A.\ref{lemma2}. To
obtain the required assertion, it remains to recall that
$H^\alpha_{n(k)}$ is malnormal in $L^\alpha_{n(k)}$ by~(v) and
take
$L=L^\alpha_{n(k)}$, $K=H^\alpha_{n(k)+1}$, and
$H=H^\alpha_{n(k)}$ in Lemma~A.\ref{lemma2}.

The topologizability of infinitely many groups $L_n^\alpha$
for every nonzero $\alpha$ implies that any
finite subset of $M$ is contained in a topologizable subgroup.
Indeed, any such subset $F$ is contained in $M_{\alpha+1}$ for some
$\alpha$. On the other hand, $M_{\alpha+1}$ is the union of the
increasing sequence of the groups $L^\alpha_n$; hence
$F$ is contained in $L_k^{\alpha}$ for some $k\in\omega$. Any
topologizable group $L_{n(k)}^\alpha$ with $n(k)\ge k$ contains $F$.

Since $M_1$ is an arbitrary non--finitely generated countable group,
any at most countable group can be embedded as a subgroup in a group
having the same properties as $M$. We obtain the following corollary.

\begin{corollary}
Any at most countable group can be embedded as a subgroup in a group
$G$ with the following properties\textup:
\begin{enumerate}
\item $G$ is an uncountable group\textup;
\item $G=\bigcup_{\alpha\in \omega_1} G_\alpha$,
where each $G_\alpha$ is a countable subgroup of
$G_{\beta}$ for any $\beta\ni\alpha$,
each $G_\alpha$ is malnormal in $G$, and all of the
$G_\alpha$ \textup(except possibly $G_1$\textup)
are topologizable\textup;
\item under CH, $G= S^{10000}$ for any uncountable $S\subset G$
\textup(this means that $G$ is a Jonsson semigroup, i.e., all proper
subsemigroups of $G$ are countable\textup) and $G$ is nontopologizable\textup;
\item under CH, $G$ is simple\textup;
\item under CH, $G\setminus\{1\}$ is
unconditionally closed but not algebraic.
\end{enumerate}
\end{corollary}

\begin{remark}
Lemma~\ref{lemma1} may
be useful for constructing an example in ZFC. The
nontopologizable group constructed by G.~Hesse in ZFC in his
dissertation~\cite{Hesse} is very likely to have such a
structure.
\end{remark}

\section*{Appendix}

We begin this section with mentioning some basic definitions
and facts from~\cite{Schupp}; see~\cite{Schupp} for more details.

\begin{extradefinition}
Suppose that $K$ and $L$ are groups, $H\subset K$ and $H'\subset L$
are their isomorphic subgroups, and $\varphi\colon H\to H'$ is an
isomorphism. The free product of $K$ and $L$ with the
subgroups $H$ and $H'$ amalgamated by the isomorphism $\varphi$ is
the quotient of the free product $K*L$ by the relations $\varphi(h) =
h$ for all $h\in H$. In what follows, we identify $H$ with $H'$
(i.e., assume that $K\cap L=H$) and refer to the
free product of $K$ and $L$ with $H$ and $H'$ amalgamated by
$\varphi$ as the \emph{free product of $K$ and $L$ with amalgamation
over $H$} or simply the \emph{amalgamated free product of $K$ and
$L$}. We use the standard notation $K\mathbin{*_H} L$ for the
amalgamated free product.
\end{extradefinition}

The groups $K$ and $L$ are naturally embedded in $K\mathbin{*_H} L$
(see~\cite{Schupp}).

We set $L^*= K\mathbin{*_H} L$ and identify the groups $K$ and $L$
with their images in $L^*$ under the natural embeddings. We refer
to elements of $L^*$ as  \emph{words} and to elements of $K$ and $L$
as \emph{letters}.

A \emph{normal form} of a nonidentity element
$w\in L^*$ is a sequence $g_1\dots g_n$ of letters such that
$w=g_1\dots g_n$ in $L^*$, $g_i$ and $g_{i+1}$ belong to different
factors ($K$ and $L$) for any $i=1, \dots, n-1$, and if $n\ne 1$,
then none of the letters $g_1,\dots, g_n$ belongs to $H$. Any element
$w$ of $L^*$ can be written in normal form.
Moreover, it may have many normal forms,
but the number of letters in each of its normal forms
is the same (see~\cite{Schupp}); it is called the
\emph{length} of $w$ and denoted by $|w|$.

\begin{extralemma}
\label{normalforms}
Any two normal forms $x_1\dots x_n$ and $y_1\dots y_n$ of the same
element of $L^*$ are related as follows: there exist $h_1, \dots,
h_{n-1}\in H$ such that $y_1=x_1h_1^{-1}$, $y_2=h_1 x_2h_2^{-1}$,
$y_3=h_2 x_3h_3^{-1}$, \dots, $y_n=h_{n-1} x_n$.
\end{extralemma}

\begin{proof}
We have $y_n^{-1}\dots y_1^{-1}x_1\dots x_n =1$. The normal form
theorem for amalgamated free products~\cite[Theorem~IV.2.6]{Schupp}
asserts that if $z_1\dots z_n$ is a normal form of some word, then
either $n=1$ and $z_1=1$ or this word is not $1$. Thus,
$y_n^{-1}\dots y_1^{-1}x_1\dots x_n =1$ is not a normal form, i.e.,
the letters $y_1^{-1}$ and $x_1$ belong to the same factor.
For definiteness, we assume that $x_1, y_1^{-1}\in
K$. Suppose that $y_1^{-1}x_1\notin H$. Let $z = y_1^{-1}
x_1$.  Since the forms $x_1\dots x_n$ and $y_1\dots y_n$ are normal,
it follows that $x_2, y_2^{-1}\in L\setminus H$.  Therefore,
$y_n^{-1}\dots y_2^{-1}zx_2\dots x_n$ is a normal form, which
contradicts its being equal to~$1$. Thus,
$y_1^{-1}x_1=h_1$ for some $h_1\in H$, whence
$y_1=x_1h_1^{-1}$. We set $y'_2=h_1^{-1}y_2$. Consider the word
$y_n^{-1}\dots {y'_2}^{-1}x_2\dots x_n$. It equals $1$; therefore, it
is not a normal form. Arguing as above, we conclude that
${y'_2}^{-1}$ and $x_2$ cancel each other, i.e.,
${y'_2}^{-1}x_2=h_2\in H$, i.e., $y_2^{-1}h_1x_2= h_2$, whence
$y_2=h_1 x_2h_2^{-1}$.  Continuing, we obtain the required $h_1,
\dots, h_n$.
\end{proof}

A word $w$ is said to be
\emph{cyclically reduced} if it has a normal form $g_1\dots g_n$ such
that $n\le 1$ or $g_1$ and $g_n$ belong to different factors
(Lemma~A.\ref{normalforms} implies that any normal form of a
cyclically reduced word has this property).  A word $w=g_1\dots g_n$
in normal form is \emph{weakly cyclically reduced} if $n\le 1$ or
$g_ng_1\notin H$.

Let $u$ and $v$ be words with normal forms $u=g_1\dots g_n$ and
$v=h_1\dots h_m$. If $g_nh_1\in H$, then we say that
$g_n$ and $h_1$ \emph{cancel} each other in
the product $uv$.  If $g_n$ and $h_1$ belong to the same factor but
$g_nh_1\notin H$, then we say that $g_n$ and $h_1$ \emph{merge} in
the normal form of the product $uv$. A representation $u_1\dots u_k$
(where the $u_i$ are words) of a word $w$ is \emph{semireduced} if
there are no cancellations in the product $u_1\dots u_k$; mergings
are allowed. If the product contains neither cancellations nor
mergings, then the representation is said to be \emph{reduced}.

A subset $R$ of the group $L^*$ is called \emph{symmetrized} if $r\in
R$ implies that $r$ is weakly cyclically reduced and all weakly
cyclically reduced conjugates of $r$ and $r^{-1}$ belong to $R$. The
\emph{symmetrized closure} of an element (or a set of elements)
of $L^*$ is the
least symmetrized set containing this element (or set).
A word $b$ is called a
\emph{piece} (with respect to a symmetrized set $R$) if there exist
different $r, r'\in R$ and some $c, c'\in L^*$
such that $r=bc$, $r=bc'$, and these representations are semireduced.

Let $\lambda > 0$.

We say that a symmetrized set $R$ satisfies the \emph{small
cancellation condition} $C'(\lambda)$ if it has the following
property.

\begin{scc}
If $r\in R$ has a semireduced representation $r=bc$, where $b$ is a
piece, then $|b|<\lambda|r|$; moreover, $|r|>1/\lambda$ for all
$r\in R$.
\end{scc}

\begin{extralemma}
\label{scc}
Suppose that $x$ and $y$ are good fellows in $K$ over $H$,
$a\in L\setminus H$, $a^{-1}H a \cap H=\{1\}$, and $h\in K\setminus
H$. Then the symmetrized closure $R$ of the  word
$$
r_0=h a ya xa (ya)^2 xa (ya)^3 \dots xa(ya)^{80}
$$
satisfies the condition $C'(\lambda)$.
\end{extralemma}

\begin{proof}
Clearly, any weakly cyclically reduced element of the group $L^*$ is
conjugate to a cyclically reduced element by means of an element of
$K\cup L$. By Theorem~IV.2.8 from \cite{Schupp}, any cyclically
reduced element of $R$ is conjugate to a cyclic permutation of
$r_0^{\pm 1}$ by  means of an element of $H$. Thus, any element of
$R$ is conjugate to a cyclic permutation of $r_0^{\pm
1}$ by means of an element of $K\cup L$ and hence has
length $6640$ ($=|r_0|$) or $6641$.

Take two elements $r, r'\in R$. Let us show that if they have
normal forms in which the initial fragments of length larger than
600 coincide, then these elements themselves coincide. Suppose that
$$
r= z_0z_1\dots z_n \qquad\text{and} \qquad r'= z'_0z'_1\dots z'_n
$$
are normal forms and $z_i=z'_i$ for $i=0,1, \dots, s$, where $s
\ge 600$. We have
$$
z_0z_1\dots z_n = t\tilde z_1 \dots \tilde z_n t^{-1}\qquad\text{and}
\qquad
z'_0z'_1\dots z'_n = t'\tilde z'_1 \dots \tilde z'_n {t'}^{-1},
$$
where $t, t'\in K\cup L$ and $\tilde z_1 \dots \tilde z_n$ and
the words $\tilde z'_1 \dots \tilde z'_n$ are cyclic permutations of
$r_0^\varepsilon$ and $r_0^\delta$ for some $\varepsilon,
\delta = \pm1$. For definiteness, suppose that $\delta=1$. Clearly,
we can assume that $t$ and $\tilde z_1$ belong to different factors
(otherwise, we replace $t$ by $t\tilde z_1$ and consider the cyclic
permutation $\tilde z_2 \dots \tilde z_n\tilde z_1$ of
$r_0^\varepsilon$); similarly, we can assume that $t$ and $\tilde
z_1$ belong to different factors as well. Then $\tilde z_n$ and
$t^{-1}$ belong to the same factor, i.e.,  $\tilde z_nt^{-1}=u\in
K\cup L$, and $t\tilde z_1 \dots \tilde z_{n-1}u$ is a
normal form. Similarly, $t'\tilde z'_1 \dots \tilde
z'_{n-1}u'$ is a normal form for some $u\in K\cup L$. By
Lemma~A.\ref{normalforms}, there exist $\tilde h_0, \dots,
\tilde h_{s}, \tilde h'_0, \dots, \tilde h'_{s}\in H$ for
which
\begin{equation}
\begin{aligned}
\label{equality}
t{\tilde h_0}^{-1}=z_0&=z'_0=t'{\tilde {h'}_0}^{-1}, \\
\tilde h_0 \tilde z_1{\tilde h_1}^{-1}=z_1&=z'_1=\tilde h'_0 \tilde
z'_1{\tilde {h'_1}}^{-1},\\
\tilde h_1 \tilde z_2{\tilde h_2}^{-1}=z_2&=z'_2=\tilde h'_1 \tilde
z'_2{\tilde {h'_2}}^{-1},\\
&\dots, \\
\tilde h_{s-1} \tilde z_{s}{\tilde h_{s}}^{-1}=
z_{s}&=z'_{s}=\tilde
h'_{s-1} \tilde z'_{s}{\tilde {h'_{s}}}^{-1}.
\end{aligned}
\end{equation}
Hence there exist $h_0, \dots,
h_{s}\in H$  such that
$$
t'=th_0^{-1}\quad\text{and}\quad
\tilde z'_i= h_{i-1} \tilde z_ih_i^{-1}\quad
\text{for}\quad i\le s.
$$
Each of the letters $z_i$ and $z'_i$ is $x^{\pm 1}$, $y^{\pm 1}$,
$a^{\pm 1}$, or $h^{\pm 1}$. Since $x$ and $y$ are good fellows over
$H$ and $x, y, h\in K\setminus H$, while $a\in L\setminus H$, it
follows that (i)~$\tilde z_i=a^\varepsilon \iff \tilde
z'_i=a$; (ii)~$\tilde z_i=x^\varepsilon$ or $\tilde z_i=h^\varepsilon \iff
\tilde z'_i=x$ or $\tilde z'_i= h$; (iii)~$\tilde z_i=y^\varepsilon$ or
$\tilde z_i=h^\varepsilon \iff \tilde z'_i=y$ or $\tilde z'_i=h$.

Suppose that $\tilde z_i\ne \tilde z'_i$, i.e., (ii) or (iii) holds.
For definiteness, we assume that $\tilde z_i=x^\varepsilon$ and
$\tilde z'_i=h$. If $i\le s-8$, then $\tilde
z'_{i+2} = y$, $\tilde z'_{i+4} = x$,
$\tilde z'_{i+6} = y$, and $\tilde z'_{i+8} = y$,
while, certainly, either $\tilde z_{i+2} = \tilde z_{i+4} = y^\varepsilon$
or $\tilde z_{i+2} = y^\varepsilon$, $\tilde z_{i+4} =
h^\varepsilon$, $\tilde z_{i+6} = y^\varepsilon$, and $\tilde z_{i+8}
= x^\varepsilon$. If $i>s-8$, then $\tilde z'_{i-2}= \tilde z'_{i-4}
= \dots = \tilde z'_{i-160}= y$, while at least one of the
corresponding letters $\tilde z_j$ is $x^\varepsilon$. In any case,
there exists a $j\le s$  such that $\tilde z_j= x^\varepsilon$ and
$\tilde z'_j=y$ or $\tilde z_j= y^\varepsilon$ and $\tilde z'_j=x$,
which is impossible.

Thus, we have $\tilde z_i^\varepsilon=\tilde z'_i$ for any $i\le
s$. Clearly, the word $\tilde z'_1 \dots \tilde z'_{s}$
(being a cyclic permutation of $r_0$) contains a fragment of the form
$xa(ya)^kxa(ya)^{k+1}xa$. The corresponding fragment of the word
$\tilde z_1\dots \tilde z_{s}$ must have the form
$x^\varepsilon a^\varepsilon (y^\varepsilon a^\varepsilon
)^kx^\varepsilon  a^\varepsilon (y^\varepsilon a^\varepsilon
)^{k+1}x^\varepsilon a^\varepsilon $, which implies $\varepsilon =
1$. These fragments, together with their positions in the words
$\tilde z_1\dots \tilde z_{s}$ and $\tilde z'_1\dots \tilde
z'_{s}$ (which are initial fragments of cyclic permutations of
$r_0$), uniquely determine the permutations. We conclude that $\tilde
z_1\dots \tilde z_{s}$ coincides with $\tilde z'_1\dots \tilde
z'_{s}$. It remains to show that $t=t'$.

As mentioned above (see~\eqref{equality}, there exist $\tilde h_0,
\tilde h_1, \tilde h_2, \tilde h'_0, \tilde h'_1, \tilde h'_2\in H$
such that
\begin{align*}
t{\tilde h_0}^{-1}&=t'{\tilde {h'}_0}^{-1}, \\
\tilde h_0 \tilde z_1{\tilde h_1}^{-1}&=\tilde h'_0 \tilde
z'_1{\tilde {h'_1}}^{-1}\quad
\text{(i.e., $\tilde z_1^{-1}\tilde h_0^{-1}\tilde h'_0 \tilde
z'_1\in H$)},\\
\intertext{and}
\tilde h_1 \tilde z_2{\tilde h_2}^{-1}&=\tilde h'_1 \tilde
z'_2{\tilde {h'_2}}^{-1}\quad
\text{(i.e., $\tilde z_2^{-1}\tilde h_1^{-1}\tilde h'_1 \tilde
z'_2\in H$)}.
\end{align*}
One of the letters $\tilde z_1=\tilde z'_1$ and $\tilde z_2=\tilde
z'_2$ is $a$. If $\tilde z_1=\tilde z'_1=a$, then $h_0=h'_0$ (because
$a^{-1}Ha\cap H=\{1\}$ by assumption), whence $t=t'$; if
$\tilde z_2=\tilde z'_2=a$, then $h_1=h'_1$, whence $h_0=h'_0$
(because $\tilde z_1=\tilde z'_1$) and $t=t'$.

Let $b$ be a piece. This means by definition
such that $b$ has two normal
forms coinciding (up to their last letters) with initial fragments
of normal forms of two different element $r$ and $r'$ in $R$; i.e.,
that there are different normal forms $z_0z_1\dots z_n$
and $z'_0z'_1\dots z'_n$ in $R$  such that
$b= z_0z_1\dots z_su=z'_0z'_1\dots z'_su'$, where $s< n$ and $u$
and $u'$ are some (possibly identity) letters.
We have shown that $s<600$ (otherwise, the forms $z_0z_1\dots z_n$
and $z'_0z'_1\dots z'_n$ would coincide).
It follows that $|b|\le 601<\frac1{10}6640$.
It remains to recall that all elements of $R$ have length 6640 or
6641.
\end{proof}

Theorem V.11.2 from \cite{Schupp} asserts, in particular, that if $N$
is the normal closure of a symmetrized set $R$ in $L^*=
K\mathbin{*_H}L$ and $R$ satisfies the condition $C'(1/10)$, then the
natural homomorphism $L^*\to L^*/N$ acts as an
endomorphism on $K$ and $L$; moreover, any nonidentity
element $w$ of $N$ has a reduced representation $w=usv$,
where $|s|>\frac 7{10}|r|$
for some $r\in R$ (and hence $|w|>7$) and $r$ has a reduced
representation of the form $r=st$.

Let $\varphi\colon L^*\to L^*/N$ be the natural homomorphism.

\begin{extralemma}\label{malnormal}
If the conditions of Lemma~\textup{A.\ref{scc}} hold and $H$ is malnormal in
$L$, then $\varphi(K)$ is
malnormal in $L^*/N$.
\end{extralemma}

\begin{proof}
Suppose that $\varphi(K)$ is not malnormal in $L^*/N$.
Take $u\in L^*$ such that $\varphi(u)\in
L^*/N\setminus \varphi(K)$ (i.e., $u\notin KN$) and
$\varphi(u)^{-1}\varphi(g)\varphi(u)=\varphi(g')$
for some $g, g'\in K\setminus \{1\}$.  This means that
$u^{-1}gu{g'}^{-1}\in N$ for some $g, g'\in K\setminus \{1\}$, or,
equivalently, $gu^{-1}g'ug''\in N$ for some $g, g', g''\in
K$ such that $g'\ne 1$ and $gg''\ne 1$. Suppose that $u$ is a
shortest word from $L^*\setminus KN$ for which such $g$, $g'$, and
$g''$ exist.  Let $u_1\dots u_n$ be a normal form of $u$. If
$u_n\in K$, then
$gu_n^{-1}=gg''{g''}^{-1}u_n^{-1}\in K$ and $u_ng''=u_ng^{-1}gg''\in
K$; replacing $g''$ by $u_ng''$ and $g$ by $gu_n^{-1}$, we see
that $u_1\dots u_{n-1}$ is a word with the same properties as $u$ but
shorter than $u$. Thus, $u_n\notin K$, i.e.,
$u_n\in L\setminus H$.

If $u_1^{-1}g'u_1=1$, then $gu^{-1}g'ug''=gg''$. As mentioned above, any
nonidentity element of $N$ has length at least~7; hence $gg''=1$,
which contradicts the assumption. Therefore, $u_1^{-1}g'u_1\ne 1$.
If $u_1\in K$, then, replacing $g'$ by $u_1^{-1}g'u_1$,
we see that $u_2\dots u_{n}$ is a word with the
same properties as $u$ but shorter than $u$. Thus, $u_1\notin
K$, i.e., $u_1\in L\setminus H$.

If $u$ has a
reduced representation $vsw$, where $s$ is a fragment of some $r\in
R$ (i.e., $r$ has a reduced representation $r=s_1ss_2$), then
$\varphi(u)=\varphi(vs_1^{-1}s_2^{-1}w)$, because
$vs_1^{-1}s_2^{-1}s^{-1}v^{-1}\in N$ (the element
$s_1^{-1}s_2^{-1}s^{-1}$ is a cyclic permutation of $r^{-1}=
s_2^{-1}s^{-1}s_1^{-1}$ and hence belongs to $R$). Thus, we have
$|s|\le |s_1|+|s_2|$ (otherwise, the word $u$ is not shortest);
i.e., $u$ cannot contain a fragment of a word $r\in R$ of length
$>\frac12|r|$.

Let us find a normal form of $gu^{-1}g'ug''$. If $g, g', g''\notin
H$ (i.e., $g, g', g''\in K\setminus H$),
then $gu_n^{-1}\dots u_1^{-1}g'u_1\dots u_ng''$ is a normal form,
because, as shown above, $u_1, u_n\in L\setminus H$.
If $g\in H$ and $g', g''\notin H$, then ${u'}_n^{-1}\dots
u_1^{-1}g'u_1\dots u_ng''$, where ${u'}_n^{-1}=gu_n^{-1}$, is a
normal form (clearly, $gu_n^{-1}\in L\setminus H$). If $g'\in H$ and
$g, g''\notin H$, then $gu_n^{-1}\dots u_2^{-1}u_0u_2\dots u_ng''$,
where $u_0=u_1^{-1}g'u_1$, is a normal form. (Indeed, we have $u_1\in
L\setminus H$ and $g'\in H\setminus \{1\}$; since $H$ is malnormal in
$L$, it follows that $u_1^{-1}g'u_1\in L\setminus H$.) The remaining
cases are considered similarly.

Thus, in any case, $gu^{-1}g'ug''$ has a normal form equal (up to
the first and last letters) to $u_n^{-1}\dots
u_2^{-1}\tilde uu_2\dots u_n$, where $\tilde u$ is the word
$u_1^{-1}g'u_1$ or the letter from $L\setminus H$ equal to
$u_1^{-1}g'u_1$.

As mentioned above, any nonidentity element of $N$ is a reduced
product of a fragment $s$ of some word $r\in R$ of length
$>\frac7{10}|r|$ and something else. Every $r\in R$ is a cyclic
permutation of $r_0$ conjugate by means of some letter. Thus,
the normal form of $gu^{-1}g'ug''$ contains a long (of length
$>\frac7{10}|r_0|-2$) fragment $t$ of a cyclic permutation of $r_0$.
Since $u$ can contain only fragments of length $\le \frac12|r|$, it
follows that $t=u_k^{-1}\dots u_2^{-1}\tilde uu_2\dots u_m$, where
$k, m> \frac1{10}|r|$. Let the fragment $t$ be
$z_1\dots z_l$, where $l=k+m+1$ or $l=
k+m+3$ (depending on $\tilde u$). According to
Lemma~A.\ref{normalforms}, for each $i\le l$, the $i$th letter in
$u_k^{-1}\dots u_2^{-1}\tilde uu_2\dots u_m$ belongs to
$Hz_iH$. Since $k$ and $m$ are large and $\tilde u$
contains one or three letters, there exists a $j\in \{2, ...  ,
\min\{k,m\}\}$ such that $u_j^{-1}\in Hx^{\pm1}H$ and $u_j\in
Hy^{\pm1}H$ ($x$ and $y$ are the same as in Lemma~A.\ref{scc}). This
contradicts the $x$ and $y$ being good fellows over $H$.
\end{proof}

We identify $K$ with $\varphi(K)$ and $L$ with $\varphi(L)$, that is,
treat $K$ and $L$ as subgroups of $(L\mathbin{*_H} K)/N$.

The following fact was
kindly communicated to the author by Anton Klyachko.

\begin{extralemma}\label{lemma2}
Suppose that $L$ and $K$ are infinite countable groups, $L\cap K =H$,
$x,y\in K$ are good fellows over $H$, $a\in L$,
$a^{-1}Ha\cap H=\{1\}$ in $L$,
$h\in K\setminus H$,
and
$$
r_0=hayaxa(ya)^2xa(ya)^3
\dots xa(ya)^{80}\in L\mathbin{*_H} K.
$$
Let $R$ be the symmetrized closure of $r_0$, and let $N$ be the
normal closure of $R$. Then the group $\langle L\mathbin{*_H} K \mid
r_0=1\rangle = (L\mathbin{*_H} K)/N$ admits a nondiscrete Hausdorff group
topology.
\end{extralemma}

\begin{proof}
Let us enumerate the elements of $L\mathbin{*_H} K$:
$$
L\mathbin{*_H} K=\{1, g_1, g_2, \dots\}.
$$
We shall construct nontrivial normal subgroups $N_1$, $N_2$,
\dots of $L\mathbin{*_H} K$ such that $N_{i+1}\subset N_i$ and $g_i\notin N_i$
for each $i$.

Take cyclically
reduced words $r_n$ in $L\mathbin{*_H}K$ such that their lengths
unboundedly increase and the symmetrized closure of $\{r_n:n\ge
0\}$ (and, therefore, of any subset of this set)
satisfies $C'(1/10)$; in particular,
each word in the normal subgroup generated by the
(symmetrized closure of) $\{r_n:n\ge k\}$
is at least half as long as $r_k$.
For such words we can take
$$
r_n=xa(ya)^{80(n-1)+1}xa(ya)^{80(n-1)+2}\dots xa(ya)^{80n}.
$$
This is proved in precisely the same way as Lemma~A.\ref{scc}.
The only difference is that $\tilde z_1 \dots \tilde z_n$ and
$\tilde z'_1 \dots \tilde z'_n$ may be cyclic permutations of
$r_k^\varepsilon$ and $r_{k'}^\delta$ for different $k$ and $k'$.
This does not matter, because if, say,
$|r_k^\varepsilon|\le |r_{k'}^\delta|$
 and $|s|> \frac1{10} |r_k^\varepsilon|$,
then the word $\tilde z_1 \dots \tilde z_s$
(as well as $\tilde z'_1 \dots \tilde z'_s$), being a
cyclic permutation of $r_k^\varepsilon$, still
contains a fragment of the form
$$
xa(ya)^jxa(ya)^{j+1}xa\quad \text{or}\quad
a^{-1}x^{-1}(a^{-1}y^{-1})^{j+1}a^{-1}x^{-1}(a^{-1}y^{-1})^ja^{-1}x^{-1},
$$
which determines $k$, $\varepsilon$, and the permutation.

For every $n\in \omega$, let $k(n)$
be an integer such that the word $r_{k(n)}$ is twice as long as $g_n$; we
assume that $k(n+1)>k(n)$.  We define $N_n$ to be the normal subgroup
generated by $\{r_k: k\ge k(n)\}$. It
does not contain $g_n$, because $g_n$ is too short.  Therefore,
$\bigcap N_n=\{1\}$. On the other hand, $N_n\not\subset NN_{n+1}$ for any
$n$; for example, $r_{k(n)}\notin NN_{n+1}$ for any $n>0$. Indeed,
$NN_{n+1}$ is the normal closure of
the set $\{r_0\}\cup\{r_k:k\ge k(n+1)\}$ and, therefore, of
the symmetrized closure $R_{n+1}$ of this set, which satisfies
the condition $C'(1/10)$.
By above-cited Theorem~V.11.2 from \cite{Schupp},
each element of $NN_{n+1}$ must contain a fragment $s$
of some $r\in R_{n+1}$
of length at least $\frac7{10}|r|$, while $r_{k(n)}$
cannot contain such a fragment.
Indeed, if $r_{k(n)}= usv$ is a reduced representation and $s$ is a long
fragment of $r$, i.e., $r$ has a reduced representation
$u'sv'$, then $svu$ (which is a cyclic permutation of $r_{k(n)}$)
is a reduced representation of some word
$\tilde r$ from the symmetrized closure of $r_{k(n)}$, and $sv'u'$
(which is a cyclic permutation of a weakly cyclically reduced conjugate
of some word in $\{r_0\}\cup\{r_k:k\ge k(n+1)\}$)
is a semireduced representation of some word ${\tilde r}'$ from
the symmetrized closure $R_{n+1}$
of $\{r_0\}\cup\{r_k:k\ge k(n+1)\}$. Thus,
$s$ is a piece
with respect to the symmetrized closure of $\{r_n:n\ge 0\}$
(which satisfies the small cancellation condition
$C'(1/10)$), and it cannot be longer than $\frac1{10}|{\tilde r}'|$.
Clearly, $|{\tilde r}'|\le |r|+1$, and
$|s|\le\frac1{10}|{\tilde r}'|<\frac 7{10}|r|$, which contradicts the
choice of $s$.

Thus, the images of the groups $N_n$ under the natural homomorphism
$L\mathbin{*_H} K\to (L\mathbin{*_H} K)/N$ form a strictly
decreasing sequence of nontrivial normal
subgroups with trivial intersection. Clearly,
such subgroups constitute a neighborhood base at the identity for some
nondiscrete Hausdorff group topology on $(L\mathbin{*_H} K)/N$.
\end{proof}

\section*{Acknowledgments}

The author thanks Anton~A.~Klyachko for many helpful discussions and
advice and the referee for very useful comments.


\begin{thebibliography}{0}

\bibitem{Markov1945}
A. A. Markov, ``On free topological groups,'' Izv. Akad. Nauk SSSR,
Ser. Mat. \textbf{9} (1), 3--64 (1945); English translation:
``Three papers on topological groups:
I. On the existence of periodic  connected topological groups.
II. On free topological groups.
III. On  unconditionally closed sets,''
Amer. Math. Soc. Transl. \textbf{30} (1950).

\bibitem{Markov1946}
A. A. Markov, ``On unconditionally closed sets,'' Mat. Sb.
\textbf{18(60)} (1), 3--26 (1946); English translation:
``Three papers on topological groups:
I. On the existence of periodic  connected topological groups.
II. On free topological groups.
III. On  unconditionally closed sets,''
Amer. Math. Soc. Transl. \textbf{30} (1950).

\bibitem{Markov1944}
A. A. Markov, ``On unconditionally closed sets,'' Dokl. Akad. Nauk SSSR,
\textbf{44} (5), 196--197 (1944); English translation:
A. Markoff,
``On unconditionally closed sets,''
C. R. (Doklady) Acad. Sci. URSS (N.S.)
\textbf{44}, 180--181  (1944).

\bibitem{new}
O. V. Sipacheva, ``Unconditionally closed and algebraic sets in
subgroups of direct products of countable groups,''
ArXiv:math.GR/0610430.

\bibitem{Shelah}
S. Shelah, ``On a problem of Kurosh, Jonsson groups, and
applications,'' in \textit{Word Problems II}, Ed. by S.~I.~Adian,
W.~W.~Boone, and G.~Higman (North-Holland, Amsterdam, 1980),
pp.~373--394.

\bibitem{Hesse}
G. Hesse,  \textit{Zur Topologisierbarkeit von Gruppen},
Dissertation (Univ. Hannover, Hannover, 1979).

\bibitem{Ol'shanskii}
A. Yu. Ol'shanskii,
``A remark on a countable nontopologized group,'' Vestnik Moskov. Univ.
Ser.~I Mat. Mekh., No.~3, p.~103 (1980);
\textit{The Geometry of  Defining Relations in Groups}
(Nauka, Moscow, 1989; Kl\"uwer, Dordrecht, 1991).

\bibitem{nontop}
Anton A. Klyachko and Anton V. Trofimov,
``The number of non-solutions of an equation in a group,''
J. Group Theory \textbf{8} (6), pp.~747--754 (2005).

\bibitem{Schupp}
R. C. Lyndon and P. E. Schupp,
\textit{Combinatorial Group Theory}
(Springer-Verlag, Berlin--Heidelberg--New York,
1977).

\end{thebibliography}
\end{document}